\documentclass[12pt,a4paper,english,reqno]{amsart}
\usepackage[a4paper,footskip = 1.5em]{geometry}
\usepackage{amsmath,amssymb,amsthm,mathtools,mathrsfs}
\usepackage{adjustbox,tikz,calc,graphics,babel}
\usepackage{csquotes,enumerate,verbatim}
\usepackage[final]{microtype}
\usepackage[numbers]{natbib}
\usetikzlibrary{shapes.misc,calc,intersections,patterns,decorations.pathreplacing}
\usepackage{hyperref}
\hypersetup{colorlinks=true, linkcolor=blue, citecolor=blue}

\theoremstyle{plain}
\newtheorem*{theorem*}{Theorem}
\newtheorem{theorem}{Theorem}[section]
\newtheorem{lemma}[theorem]{Lemma}

\newtheorem{proposition}[theorem]{Proposition}
\newtheorem*{claim*}{Claim}
\newtheorem{corollary}[theorem]{Corollary}

\theoremstyle{remark}

\def\N{\mathbb{N}}

\def\C{\mathcal}

\let\emptyset\varnothing

\let\originalleft\left
\let\originalright\right
\renewcommand{\left}{\mathopen{}\mathclose\bgroup\originalleft}
\renewcommand{\right}{\aftergroup\egroup\originalright}

\makeatletter
\def\imod#1{\allowbreak\mkern10mu({\operator@font mod}\,\,#1)}
\makeatother

\begin{document}

\title{A canonical Ramsey theorem for exactly $m$-coloured complete subgraphs}

\author{Teeradej Kittipassorn}
\address{Department of Mathematical Sciences, University of Memphis, Memphis TN 38152, USA}
\email{t.kittipassorn@memphis.edu}

\author{Bhargav Narayanan}
\address{Department of Pure Mathematics and Mathematical Statistics, University of Cambridge, Wilberforce Road, Cambridge CB3\thinspace0WB, UK}
\email{b.p.narayanan@dpmms.cam.ac.uk}

\date{14 May 2013}

\begin{abstract}
Given an edge colouring of a graph with a set of $m$ colours, we say that the graph is (\emph{exactly}) \emph{$m$-coloured} if each of the colours is used.   We consider edge colourings of the complete graph on $\N$ with infinitely many colours and show that either one can find an $m$-coloured complete subgraph for every natural number $m$ or there exists an infinite subset $X \subset \N$ coloured in one of two canonical ways: either the colouring is injective on $X$ or there exists a distinguished vertex $v$ in $X$ such that $X \setminus \lbrace v \rbrace$ is $1$-coloured and each edge between $v$ and $X \setminus \lbrace v \rbrace$ has a distinct colour (all different to the colour used on $X \setminus \lbrace v \rbrace$). This answers a question posed by Stacey and Weidl in 1999. The techniques that we develop also enable us to resolve some further questions about finding $m$-coloured complete subgraphs in colourings with finitely many colours.
\end{abstract}

\maketitle

\section{Introduction}
\label{canon13-intro}
A classical result of Ramsey~\citep{Ramsey1930} says that when the edges of a complete graph on a countably infinite vertex set are finitely coloured, one can always find a complete infinite subgraph all of whose edges have the same colour.

Ramsey's Theorem has since been generalised in many ways; most of these generalisations are concerned with finding monochromatic substructures in various coloured structures. For a survey of many of these generalisations, see the book of Graham, Rothschild and Spencer~\citep{Graham1990}. Ramsey theory has witnessed many developments over the last fifty years and continues to be an area of active research today; see, for example,~\citep{Gyarfas2002, Leader2006, Vuksanovic2006, Larson2008}.

While one is usually concerned with finding monochromatic substructures in various finitely coloured structures, two alternative directions are as follows. First, one could study colourings that use infinitely many colours, as was first done by Erd\H{o}s and Rado~\citep{Erdos1950} and by many others after them. Second, one could look for structures which are coloured with exactly $m$ colours for some $m\ge2$. This was first considered by Erickson~\citep{Erickson1994} and then investigated further by Stacey and Weidl~\citep{Stacey1999}.  In this paper, we shall consider the question of finding structures coloured with exactly $m$ colours in colourings that use infinitely many colours.

The rest of this paper is organised as follows. In Section~\ref{canon13-results}, we present the relevant definitions that we require and the statements of our results. Section~\ref{canon13-proof} is devoted to the proof of our main result. In Section~\ref{canon13-extensions}, we describe an extension of our main result and some applications of this extension. We conclude the paper in Section~\ref{canon13-conclude} with some open problems.

\section{Our results}
\label{canon13-results}

For a set $X$, denote by $X^{(2)}$ the set of all unordered pairs of elements of $X$; equivalently, $X^{(2)}$ is the complete graph on the vertex set $X$. As usual, we write $[n]$ to denote $\{1,\dots,n\}$, the set of the first $n$ natural numbers. We denote a surjective map $f$ from a set $X$ to another set $Y$ by $f: X \twoheadrightarrow Y$. By a \emph{colouring} of a graph, we mean a colouring of the edges of the graph unless we specify otherwise.

Let $\Delta:\N^{(2)}\twoheadrightarrow \C{C}$ be a surjective colouring of the edges of the complete graph on $\N$ with an arbitrary set of colours $\C{C}$. If the set of colours $\C{C}$ is infinite, we say that $\Delta$ is an \emph{infinite-colouring} and if $\C{C}$ is finite, we say that $\Delta$ is a \emph{$k$-colouring} if $|\C{C}|=k$. We say that a subset $X$ of $\N$ is (\emph{exactly}) $m$-\emph{coloured} if $\Delta(X^{(2)})$, the set of values attained by $\Delta$ on the edges with both endpoints in $X$, has size exactly $m$. We write $\gamma_{\Delta}(X)$, or $\gamma(X)$ in short, for the size of the set $\Delta( X^{(2)} )$; in other words, every set $X$ is $\gamma(X)$-coloured.

Our main aim in this paper is to establish a canonical Ramsey theory for $m$-coloured graphs. Canonical Ramsey theory, which originates in a classical paper of Erd\H{o}s and Rado~\citep{Erdos1950}, provides results about colourings which use an arbitrary set of colours. We will need a basic canonical Ramsey theorem proved by Erd\H{o}s and Rado. To state this result, it will be convenient to introduce some notation. We say that $X\subset\N$ is \emph{rainbow coloured} if no two edges with both endpoints in $X$ receive the same colour. Also, we say that $X\subset\N$ is \emph{left coloured} if for $i,j,k,l\in X$ with $i<j$ and $k<l$, $\Delta(ij)=\Delta(kl)$ if and only if $i=k$, and the definition of \emph{right coloured} is analogous; if $X$ is left or right coloured, we say, in short, that $X$ is \emph{lexically coloured.} With these definitions in place, we can now state the canonical Ramsey theorem of Erd\H{o}s and Rado~\citep{Erdos1950}.

\begin{theorem}
\label{canon13-erdos}
For any colouring $\Delta:\N^{(2)}\twoheadrightarrow \C{C}$, there exists an infinite subset $X$ of $\N$ such that either
\begin{enumerate}
\item $X$ is $1$-coloured, or 
\item $X$ is rainbow coloured, or 
\item $X$ is lexically coloured. \qed 
\end{enumerate}
\end{theorem}

For a colouring $\Delta:\N^{(2)}\twoheadrightarrow\C{C}$ of the complete graph on $\N$ with an arbitrary set of colours, we define the set 
\[
\C{G}_{\Delta} = \left\{ \gamma_\Delta(X):X\subset \N \right\} .
\] 
Stacey and Weidl~\cite{Stacey1999} considered the following question: which natural numbers $m$ are guaranteed to be elements of $\C{G}_{\Delta}$ for every infinite-colouring $\Delta$? By considering a rainbow colouring $\Delta$ of $\N$, we see that unless $m=\binom{n}{2}$ for some $n\ge2$, $m$ is not guaranteed to be a member of $\C{G}_{\Delta}$. In the other direction, since an edge is a $1$-coloured complete graph, $\binom{2}{2}=1$ is always an element of $\C{G}_{\Delta}$. Stacey and Weidl were able to show that $\binom{3}{2}=3$ is also always an element of $\C{G}_{\Delta}$ for every infinite-colouring $\Delta$. But for $n\ge4$, they were unable to decide whether or not there exists an infinite-colouring $\Delta$ such that $\binom{n}{2}\notin\C{G}_{\Delta}$. In particular, they asked if all natural numbers of the form $\binom{n}{2}$ must be contained in $\C{G}_{\Delta}$ for every infinite-colouring $\Delta$.

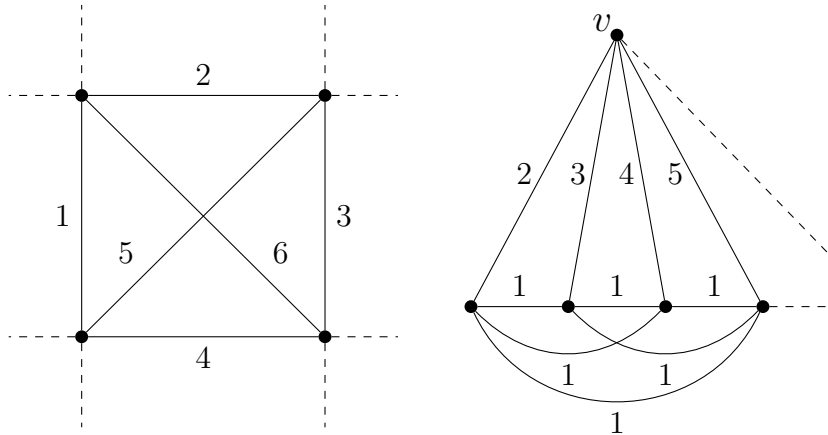
\begin{figure}
\begin{center}
\begin{tikzpicture}[auto, yscale = 0.8, xscale = 0.64]

\node (v1) at (0,0) [inner sep=0.5mm, thick, circle, draw=black!100, fill=black!100] {};
\node (v2) at (0,4) [inner sep=0.5mm, thick, circle, draw=black!100, fill=black!100] {};
\node (v3) at (5,4) [inner sep=0.5mm, thick, circle, draw=black!100, fill=black!100] {};
\node (v4) at (5,0) [inner sep=0.5mm, thick, circle, draw=black!100, fill=black!100] {};
\draw (v1) to node {1} (v2);
\draw (v2) to node {2} (v3);
\draw (v3) to node {3} (v4);
\draw (v4) to node {4} (v1);
\draw (v1) to node [near start]{5} (v3);
\draw (v2) to node [near end]{6} (v4);
\draw [dashed] (v1)--(-1.5,0);
\draw [dashed] (v1)--(0,-1.5);
\draw [dashed] (v2)--(0,5.5);
\draw [dashed] (v2)--(-1.5,4);
\draw [dashed] (v3)--(5,5.5);
\draw [dashed] (v3)--(6.5,4);
\draw [dashed] (v4)--(5,-1.5);
\draw [dashed] (v4)--(6.5,0);

\node (u1) at (8,0.5) [inner sep=0.5mm, thick, circle, draw=black!100, fill=black!100] {};
\node (c2) at (9.1,2.7) {2};
\node (u2) at (10,0.5) [inner sep=0.5mm, thick, circle, draw=black!100, fill=black!100] {};
\node (c3) at (10.2,2.7) {3};
\node (u) at (11,5) [inner sep=0.5mm, thick, circle, draw=black!100, fill=black!100] {};
\node (uu) at (10.7,5.2) {\large $v$};
\node (u3) at (12,0.5) [inner sep=0.5mm, thick, circle, draw=black!100, fill=black!100] {};
\node (c4) at (11.2,2.7) {4};
\node (u4) at (14,0.5) [inner sep=0.5mm, thick, circle, draw=black!100, fill=black!100] {};
\node (c5) at (12.2,2.7) {5};
\draw (u1) to node {1} (u2);
\draw (u1) to[out=-40,in=-140] node [swap] {1} (u3);
\draw (u2) to node {1} (u3);
\draw (u2) to[out=-40,in=-140] node [swap] {1} (u4);
\draw (u3) to node {1} (u4);
\draw (u1) to[out=-60,in=-120] node [swap] {1} (u4);
\draw (u) -- (u1);
\draw (u) -- (u2);
\draw (u) -- (u3);
\draw (u) -- (u4);
\draw [dashed] (u4)--(15.5,0.5);
\draw [dashed] (u)--(15.5,1.25);
\end{tikzpicture}
\end{center}
\caption{A rainbow colouring and a star colouring with centre $v$.}
\label{canon13-Colouring}
\end{figure}

Here, we shall consider a more general question: when is $\C{G}_{\Delta} \neq \N$?  As remarked above, for an injective colouring $\Delta$, $\C{G}_{\Delta}=\{\binom{n}{2} : n\ge2\} \neq \N$. There is another infinite-colouring $\Delta$ for which $\C{G}_{\Delta} \neq \N$ which is slightly less obvious. Given $X\subset\N$, if there is a vertex $v \in X$ such that $X \setminus \{ v \}$ is $1$-coloured and all the edges between $v$ and $X \setminus \{ v \}$ have distinct colours (which are also all different from the colour appearing in $X \setminus \{ v \}$), then we say that  $X$ is \emph{star coloured} (\emph{with centre v}). It is easy to check (see Figure~\ref{canon13-Colouring}) that if $\N$ is star coloured by $\Delta$, then $\C{G}_{\Delta}=\N\setminus\{2\}$.

Our main result, stated below, is that the two colourings described above are, in a sense, the `canonical' colourings for which $\C{G}_{\Delta}\neq\N$.

\begin{theorem}
\label{canon13-canresult}
For every infinite-colouring $\Delta:\N^{(2)}\twoheadrightarrow\N$, either 
\begin{enumerate}
\item $\C{G}_{\Delta} = \N$, or
\item there exists an infinite rainbow coloured subset of $\N$, or 
\item there exists an infinite star coloured subset of $\N$.
\end{enumerate}
\end{theorem}

An immediate consequence of Theorem~\ref{canon13-canresult} is that the answer to the question posed by Stacey and Weidl is in the affirmative.

\begin{corollary}
For every infinite-colouring $\Delta:\N^{(2)}\twoheadrightarrow\N$, and for every natural number $n\ge2$, $\binom{n}{2}\in\C{G}_{\Delta}$.\qed
\end{corollary}

We do not prove Theorem~\ref{canon13-canresult} as stated. Instead, it will be more convenient to prove a stronger result which we shall state and prove in Section~\ref{canon13-proof}.

In the context of Ramsey theory, one is usually interested in finding `large' homogeneous structures with certain properties. With this in mind, for a colouring $\Delta:\N^{(2)}\twoheadrightarrow\C{C}$, we define 
\[
\C{F}_{\Delta}=\left\{ \gamma_\Delta(X):X\subset\N \mbox{ such that } X \mbox{ is infinite}\right\} .
\]
When $\Delta$ is an infinite-colouring, it might so happen that for each infinite subset $X$ of $\N$, the set $\Delta(X^{(2)})$ is infinite; consequently, it is only really meaningful to study the set $\C{F}_{\Delta}$ in the case of colourings using finitely many colours. The question of finding $m$-coloured complete infinite subgraphs, was first considered by Erickson~\citep{Erickson1994}; see also~\citep{Stacey1999, Bhargav2013}. If $\Delta:\N^{(2)}\twoheadrightarrow[k]$ is a $k$-colouring of the edges of the complete graph on the natural numbers, then clearly $k\in\C{F}_{\Delta}$ as $\Delta$ is surjective, and Ramsey's Theorem tells us that $1\in\C{F}_{\Delta}$. Erickson~\citep{Erickson1994} noted that a fairly straightforward application of Ramsey's Theorem enables one to show that $2\in\C{F}_{\Delta}$. Erickson also conjectured that if $k>m>2$, then there is colouring $\Delta$ with exactly $k$ colours such that $m\notin\C{F}_{\Delta}$. Stacey and Weidl~\citep{Stacey1999} partially resolved this conjecture, showing for every $m>2$ that there is a $k$-colouring $\Delta$ such that $m\notin\C{F}_{\Delta}$ provided $k$ is sufficiently larger than $m$.

In the light of these negative results, Stacey and Weidl~\citep{Stacey1999} asked, in the context of colourings using finitely many colours, if we can do any better if we focus our attention on $\C{G}_\Delta$ as opposed to $\C{F}_\Delta$. More precisely, they raised the following question: do there exist natural numbers $m \in \N$ with the property that for all sufficiently large $k\in \N$,  $m \in \C{G}_\Delta$ for every $k$-colouring $\Delta:\N^{(2)}\twoheadrightarrow[k]$? Observe that any such natural number $m$, assuming one exists, must be of the form $\binom{n}{2}$ or $\binom{n}{2} + 1$ for some natural number $n\ge2$. One can see this by considering the family of `small-rainbow colourings' of the complete graph on $\N$ which colour all the edges of some finite complete subgraph with distinct colours and all the remaining edges with a single colour not used in the finite (rainbow coloured) complete subgraph. On the other hand, when $m$ is of the form $\binom{n}{2}$ or $\binom{n}{2} + 1$ for some natural number $n\ge2$, we have the following positive result.

\begin{theorem}
\label{canon13-mainfvkm}
For all $n\in\N$, there exists a natural number $C=C(n)$ such that for any $k$-colouring $\Delta:\N^{(2)}\twoheadrightarrow[k]$ with $k\ge C$, both $\binom{n}{2}, \binom{n}{2}+1 \in \C{G}_\Delta$.
\end{theorem}

It turns out that the techniques used to prove Theorem~\ref{canon13-canresult} also allow us to prove a finitary version of the same theorem. In Section~\ref{canon13-extensions}, we present this finitary result and use it prove Theorem~\ref{canon13-mainfvkm} in a slightly stronger form.

\section{Proof of the main theorem}
\label{canon13-proof}

To prove Theorem~\ref{canon13-canresult}, it will be more convenient to work with general infinite graphs. By an \emph{infinite graph}, we mean a graph whose vertex set is $\N$ and which has infinitely many edges.

It will be helpful to establish a few notational conveniences. Given an infinite graph $G$ and an infinite-colouring $\Delta:G\twoheadrightarrow\N$ of the edges of $G$, for a subset $X$ of $\N$, we shall write $\gamma_{G}(X)$, or just $\gamma(X)$ when both the colouring and the graph in question are clear from the context, for the number of distinct colours attained by $\Delta$ on $G[X]$, the subgraph of $G$ induced by $X$; if $H$ is a subgraph of $G$, we write $\gamma_H$(X) for the the number of distinct colours attained by $\Delta$ on $H[X]$. For disjoint subsets $X$ and $Y$, write $\gamma(X,Y)$ for the number of distinct colours in the induced bipartite subgraph between $X$ and $Y$ in $G$. Also, for a vertex $v\in \N$, we shall write $\gamma(v)$ for $\gamma(\{v\}, \N \setminus \{v\})$, the number of distinct colours of the edges incident to $v$ in $G$.

We define the set $\C{G}_{\Delta}$ for an infinite-colouring $\Delta:G\twoheadrightarrow\N$ of an infinite graph $G$ in the obvious way by setting
\[
\C{G}_{\Delta}=\left\{ \gamma_G(X) : X\subset \N\right\}.
\]
In a general graph $G$, we say that $X$ is \emph{rainbow coloured in $G$} if $G[X]$ is a complete subgraph of $G$ which is rainbow coloured. We say that $X$ is \emph{star coloured} (\emph{with centre $v$}) \emph{in $G$} if there is a vertex $v \in X$ such that $G[X \setminus \{ v \}]$ is either an independent set or a $1$-coloured complete graph, and all the edges between $v$ and $X \setminus \{ v \}$ are present and have distinct colours, which are also all different from the colour of $G[X \setminus \{ v \}]$ in the case where $X \setminus \{ v \}$ does not induce an independent set.The following result easily implies Theorem~\ref{canon13-canresult}.

\begin{theorem}
\label{canon13-canresult0}
For every infinite-colouring $\Delta:G\twoheadrightarrow\N$ of an infinite graph $G$, either 
\begin{enumerate}
\item $\C{G}_{\Delta} = \N$, or
\item there exists an infinite rainbow coloured subset of $\N$; or 
\item there exists an infinite star coloured subset of $\N$.
\end{enumerate}
\end{theorem}

For any finite set of colours $\C{S}$, note that if we delete all the edges of an infinite graph $G$ which are coloured with a colour from $\C{S}$ by an infinite-colouring $\Delta$ of the edges of $G$, the resulting graph $H$ is infinite and the restriction of $\Delta$ to $H$ is an infinite-colouring. This makes the statement of Theorem~\ref{canon13-canresult0} more amenable to induction than that of Theorem~\ref{canon13-canresult} and motivates the stronger statement of Theorem~\ref{canon13-canresult0}.

Fix an infinite-colouring $\Delta:G\twoheadrightarrow\N$ of an infinite graph $G$ and note that if we have a partition $X = X_1 \cup X_2 \cup \dots X_n$ of a subset $X$ of $\N$, then
\[
\sum_{1\le i \le n}{\gamma(X_i)} + \sum_{1 \le i<j \le n}{\gamma(X_i, X_j)} \ge \gamma(X).
\]
Consequently, if $\gamma(X)=\infty$, then at least one of the terms on the left is infinite; we shall make use of this fact repeatedly.
 
Next, we state a technical lemma about `almost bipartite colourings' which will be useful in proving Theorem~\ref{canon13-canresult0}.
 
\begin{lemma}
\label{canon13-biplem0}
Let $G$ be an infinite graph and suppose that an infinite-colouring $\Delta:G\twoheadrightarrow\N$ of $G$ is such that
\begin{enumerate} 
\item $\gamma(v)<\infty$ for all $v\in \N$, and 
\item there is a partition of $\N=A\cup B$ such that $\gamma(A)<\infty$, $\gamma(B)<\infty$ and $\gamma(A,B)=\infty$. 
\end{enumerate}
Then for every natural number $m$, there exists a subset $X$ of $\N$ such that $X\cap A\not=\emptyset$, $X\cap B\not=\emptyset$ and $\gamma(X)=m$.
\end{lemma}

Our strategy for proving both Theorem~\ref{canon13-canresult0} and Lemma~\ref{canon13-biplem0} is to inductively construct a set $X$ for which $\gamma_G(X)=m$. To do this, we shall first delete some edges from $G$ to get a new infinite graph $H$ so that the restriction of $\Delta$ to $H$ is also an infinite-colouring. We then inductively find a set $Y$ with $\gamma_H(Y)=l$ for a suitably chosen $l<m$. Finally, we use the deleted edges in conjunction with $Y$ to obtain $X$.

We first prove Lemma~\ref{canon13-biplem0} and then show how to deduce Theorem~\ref{canon13-canresult0} from it.
\begin{proof}[Proof of Lemma~\ref{canon13-biplem0}]
Before we begin, let us note some consequences of our assumptions about the colouring $\Delta$. Since $\gamma(v)<\infty$ for all $v\in \N$ and $\gamma(A,B)=\infty$, both $A$ and $B$ must be infinite. Furthermore, observe that if $\gamma(U)=\infty$ for some $U\subset \N$, then since $\gamma(A)<\infty$ and $\gamma(B)<\infty$, both $U\cap A$ and $U \cap B$ must be infinite.

We proceed by induction on $m$. The result is trivial for $m=1$. Assuming the result for all $l<m$, we shall prove the result for $m$. 

Pick an edge $uv$ such that $u\in A$ and $v\in B$ and say that the colour of the edge is $c$. We know that $\gamma(u)<\infty$. We may assume, relabeling colours if necessary, that the colours of the edges incident to $u$ are $1,\dots,\gamma(u)$. Consider the partition
\[
\N\setminus \{u\}=U_0\cup U_1\cup\dots\cup U_{\gamma(u)},
\]
where $U_0$ is the set of vertices not adjacent to $u$ in $G$ and for $1 \le i \le \gamma(u)$, $U_i$ is the set of all vertices that are joined to $u$ by an edge of colour $i$. By considering the following three cases, we first show that we may assume that $\gamma(U_0) = \infty$.

\textbf{Case 1: $\gamma(U_i)=\infty$ for some $i\not=0$.} We begin by observing (see Figure~\ref{canon13-Case1}) that 
\[
\gamma(U_i\cap A)+\gamma(U_i\cap B)+\gamma(U_i\cap A,U_i\cap B)\ge\gamma(U_i).
\]
Since $\gamma(U_i\cap A)\le\gamma(A)<\infty$ and $\gamma(U_i\cap B)\le\gamma(B)<\infty$, we conclude that $\gamma(U_i\cap A, U_i\cap B)=\infty$. 

Let $H$ be the infinite subgraph of $G[U_i]$ obtained by deleting all the edges of $G[U_i]$ of colour $i$. Then there exists, by the induction hypothesis, a subset $Y$ of $U_i$ such that $Y\cap (U_i\cap A)\not=\emptyset$, $Y\cap (U_i\cap B)\not=\emptyset$ and $\gamma_{H}(Y)=m-1$. Observe that all the edges between $u$ and $Y\subset U_i$ are coloured $i$ in $G$. Since the colour $i$ is not counted by $\gamma_{H}$, we see that $\gamma_{G}(Y\cup\{u\})=m$. Therefore, $X=Y\cup\{u\}$ is the required subset since $X\cap A\not=\emptyset$ and $X\cap B\not=\emptyset$.

\begin{figure}
\begin{center}
\begin{tikzpicture}[auto]

\draw [rounded corners] (0,0)--(3,0)--(3,6)--(0,6)--cycle;
\node (A) at (1,0.9) {\Large $A$};
\draw [rounded corners] (8,0)--(11,0)--(11,6)--(8,6)--cycle;
\node (B) at (10,0.9)  {\Large $B$};

\node (u) at (1.5,5.3) [inner sep=0.5mm, thick, circle, draw=black!100, fill=black!100] {};
\node (uu) at (1.25, 5.5) {\large $u$};

\draw [rounded corners] (2.5,0.75)--(8.5,0.75)--(8.5,3.5)--(2.5,3.5)--cycle;
\node (U) at (5.5,2.1)  {\Large $U_i$};

\draw (1.5,5.3)--(2.6,3.5);
\draw (1.5,5.3)--(8.4,3.5);
\draw [dashed] (1.5,5.3)--(3.6,3.5);
\draw [dashed] (1.5,5.3)--(6.4,3.5);
\node (i) at (3.6, 4.1) {\large $i$};
\draw [dashed] (u) to (i);
\draw [dashed] (i) to (4.65,3.5);

\end{tikzpicture}
\end{center}
\caption{Case 1.}
\label{canon13-Case1}
\end{figure}

\textbf{Case 2: $\gamma(U_i,U_j)=\infty$ for some $0<i<j$.} Observe (see Figure~\ref{canon13-Case2}) that $\gamma(U_i\cap A,U_j\cap A)\le\gamma(A)<\infty$ and $\gamma(U_i\cap B,U_j\cap B)\le\gamma(B)<\infty$. So we must either have $\gamma(U_i\cap A,U_j\cap B)=\infty$ or $\gamma(U_i\cap B,U_j\cap A)=\infty$. Without loss of generality, assume that $\gamma(U_i\cap A,U_j\cap B)=\infty$. 

If $m\ge3$, we may assume that the result holds for $m-2$. Let $H$ be the infinite subgraph of $G[(U_i\cap A)\cup(U_j\cap B)]$ obtained by deleting edges of colour $i$ and $j$ from $G[(U_i\cap A)\cup(U_j\cap B)]$. Then there exists, by the induction hypothesis, a subset $Y$ of $(U_i\cap A)\cup(U_j\cap B)$ such that $Y\cap (U_i\cap A)\not=\emptyset$, $Y\cap (U_j\cap B)\not=\emptyset$ and $\gamma_{H}(Y)=m-2$. Since $Y \subset U_i \cup U_j$, all the edges between $u$ and $Y$ in $G$ are coloured either $i$ or $j$, and as $Y\cap U_i\not=\emptyset$ and $Y\cap U_j\not=\emptyset$, edges of both colours are present. Since both colours $i$ and $j$ are not counted by $\gamma_{H}$, it follows that $\gamma_{G}(Y\cup\{u\})=m$. Clearly, $Y\cap A\not=\emptyset$ and $Y\cap B\not=\emptyset$, so $X=Y\cup\{u\}$ is the required subset.

Now suppose that $m=2$. Since $\gamma(w)<\infty$ for all $w\in\N$, we can greedily find an infinite matching $M=\{a_1b_1,a_2b_2,\dots \} $ between $U_i\cap A$ and $U_j\cap B$ in $G$ such that each edge of the matching has a distinct colour. If $a_k$ and $b_l$ are not adjacent in $G$ for some $k,l \in \N$, then $X=\{u,a_k,b_l\}$ is immediately seen to be 2-coloured. So we may suppose that for each $k, l \in \N$, $a_k$ is adjacent to $b_l$ in $G$. 

Since $\gamma(\{a_1,a_2,\dots\})<\infty$, it follows from Ramsey's Theorem that there exists a subset $\{a'_1,a'_2,\dots\}$ of $\{a_1,a_2,\dots\}$ which either induces an independent set or a 1-coloured complete graph. Let $a'_k$ be matched to the vertex $b'_k$ in $M$ and let $c_k$ denote the colour of the edge $a'_kb'_k$.

If $ \{ a'_1, a'_2, \dots \}$ is an independent set in $G$, then since $\gamma(a'_1)<\infty$, there exist $s,t\in\N$ such that $a'_1b'_s$ and $a'_1b'_t$ have the same colour, say $d$. By our choice of $M$,  $c_s \neq c_t$. Hence, at least one of $c_s$ or $c_t$, say $c_s$, is not equal to $d$.  Then it is easy to check that $X=\{a'_1,a'_s,b'_s\}$ is the required subset.

If $\{a'_1,a'_2,\dots\}$ induces a complete graph of colour $d$ in $G$, we may assume (by discarding the edge $a'_1b'_1$ and relabelling the remaining vertices if necessary) that $c_1$, the colour of the edge $a'_1b'_1$, is not equal to $d$. Since $\gamma(b'_1)<\infty$, there exist $s,t\in\N$ such that $\Delta(a'_sb'_1)=\Delta(a'_tb'_1)$. If $\Delta(a'_sb'_1)=d$, then we may take $X=\{a'_1,a'_s,b'_1\}$. On the other hand, if $\Delta(a'_sb'_1)\neq d$, then $X=\{a'_s,a'_t,b'_1\}$ is the required subset. 

\begin{figure}
\begin{center}
\begin{tikzpicture}[auto]

\draw [rounded corners] (0,0)--(3,0)--(3,6)--(0,6)--cycle;
\node (A) at (1,0.9) {\Large $A$};
\draw [rounded corners] (8,0)--(11,0)--(11,6)--(8,6)--cycle;
\node (B) at (10,0.9)  {\Large $B$};

\node (fin1) at (1.25, 3) {finite};
\node (fin2) at (9.75, 3) {finite};
\node (infin) at (5.5, 3) {infinite};

\draw [rounded corners] (2.5,0.5)--(8.5,0.5)--(8.5,2)--(2.5,2)--cycle;
\node (Uj) at (5.5,1.25)  {\Large $U_j$};

\draw [rounded corners] (2.5,4)--(8.5,4)--(8.5,5.5)--(2.5,5.5)--cycle;
\node (Ui) at (5.5,4.75)  {\Large $U_i$};

\draw [dashed, out = 90, in = 180] (fin1) to (2.5,4.75);
\draw [dashed, out = -90, in = 180] (fin1) to (2.5,1.25);

\draw [dashed, out = 90, in = 0] (fin2) to (8.5,4.75);
\draw [dashed, out = -90, in = 0] (fin2) to (8.5,1.25);

\draw [dashed] (infin) to (2.75, 4.75);
\draw [dashed] (infin) to (8.25, 4.75);
\draw [dashed] (infin) to (2.75, 1.25);
\draw [dashed] (infin) to (8.25, 1.25);
\end{tikzpicture}
\end{center}
\caption{Case 2.}
\label{canon13-Case2}
\end{figure}

\textbf{Case 3: $\gamma(U_0, U_i)=\infty$ for some $i\not=0$.} We argue as we did in Case 2. We may assume that $\gamma(U_0\cap A, U_i\cap B)=\infty$. Let $H$ be the infinite subgraph of $G[(U_0\cap A)\cup(U_i\cap B)]$ obtained by deleting all the edges of colour $i$ from $G[(U_0\cap A)\cup(U_i\cap B)]$.

By the induction hypothesis, there exists a subset $Y$ of $(U_0\cap A)\cup(U_i\cap B)$ such that $Y\cap (U_0\cap A)\not=\emptyset$, $Y\cap (U_i\cap B)\not=\emptyset$ and $\gamma_{H}(Y)=m-1$. As before, every edge between $u$ and $Y$ is coloured $i$ in $G$ (and $u$ is adjacent to at least one vertex of $Y$ since $Y\cap (U_i\cap B)\not=\emptyset$). Since the colour $i$ is not counted by $\gamma_{H}$, it follows that $\gamma_{H}(Y\cup\{u\})=m$. Hence, $X=Y\cup\{u\}$ is the required subset. 

Hence, we may now assume that $\gamma(U_0) = \infty$. Since $\gamma(U_0) = \infty$, $U_0$ clearly meets both $A$ and $B$ in infinitely many vertices. We consider the graph induced by $U_0 \cup \{v\}$ and let $V_0$ be the set of those vertices of $U_0$ not adjacent to $v$ in $G[U_0 \cup \{v\}]$. Since $\gamma(v)<\infty$, we have a partition of $U_0 \setminus V_0 = V_1\cup\dots\cup V_n$, with $n \le \gamma(v)$, based on the colour of the edge joining a given vertex of $U_0 \setminus V_0$ to the vertex $v$. Applying the same argument as in Cases 1, 2 and 3 (which depended only on the vertex $u$ and not on $v$) to the vertex $v$ in $G[U_0 \cup \{v\}]$, we see that we are done unless $\gamma(V_0)=\infty$. 

In this case, we consider the partition $V_0=(V_0\cap A)\cup(V_0\cap B)$. Note that $\gamma(V_0\cap A)<\infty$, $\gamma(V_0\cap B)<\infty$ and $\gamma(V_0\cap A,V_0\cap B)=\infty$. Recall that we chose $u\in A$ and $v\in B$ such that the edge $uv$ has colour $c$. Let $H$ be the infinite subgraph of $G[V_0]$  obtained by deleting edges of colour $c$ from $G[V_0]$. By the induction hypothesis, there is a subset $Y$ of $V_0$ such that $\gamma_{H}(Y)=m-1$. Observe that $uv$ has colour $c$ and furthermore, $u$ and $v$ are not adjacent to any of the vertices of $Y$. Since the colour $c$ is not counted by $\gamma_{H}$, we see that $\gamma_{G}(Y\cup\{u,v\})=m$. Therefore, $X=Y\cup\{u,v\}$ is the required subset since clearly, $X\cap A\not=\emptyset$ and $X\cap B\not=\emptyset$. This completes the proof. 
\end{proof} 

We are now in a position to deduce Theorem~\ref{canon13-canresult0} from Lemma~\ref{canon13-biplem0}.

\begin{proof}[Proof of Theorem~\ref{canon13-canresult0}]
Let $\Delta:G\twoheadrightarrow\N$ be an infinite-colouring of an infinite graph $G$. We shall prove by induction on $m$ that if $G$ contains no infinite rainbow coloured or star coloured subset, then $m\in \C{G}_{\Delta}$ for each $m \in \N$. The result is trivial for $m=1$. Now suppose that $m\ge2$. We shall inductively find a subset $X$ of $\N$ with $\gamma(X)=m$.

If $\gamma(v)=\infty$ for some vertex $v\in \N$, then we can find an infinite subset $U=\{u_1,u_2,\dots\}$ of $\N$ such that the edges $vu_i$ and $vu_j$ have distinct colours for all $i\not=j$. Applying Theorem~\ref{canon13-erdos} to the restriction of $\Delta$ to $G[U]$ (by colouring non-edges with a new colour, for example),  we can find an infinite subset $W=\{w_1,w_2,\dots\}$ of $U$ such that $W$ is either an independent set, 1-coloured, rainbow coloured or lexically coloured. By assumption, $W$ cannot be rainbow coloured. If $W$ is either an independent set or 1-coloured, it is clear that $ W\cup\{v\}$ is star coloured with centre $v$. If $W$ is lexically coloured, then it is easy to check that $\C{G}_\Delta = \N$.

So we may assume that $\gamma(v)<\infty$ for all $v\in \N$. Pick an edge $uv$ of $G$, and say that the colour of the edge is $c$.  We may suppose that the colours of the edges incident to $u$ are $1,\dots,\gamma(u)$. Consider the partition
$\N\setminus \{u\}=U_0\cup U_1\cup\dots\cup U_{\gamma(u)}$, where $U_0$ is the set of vertices not adjacent to $u$ in $G$ and for $1 \le i \le \gamma(u)$, $U_i$ is the set of all vertices that are joined to $u$ by an edge of colour $i$. Since $\gamma(\N)=\infty$, by the pigeonhole principle, we must either have $\gamma(U_i)=\infty$ for some $i$, or $\gamma(U_i,U_j)=\infty$ for some $i\not= j$. We distinguish the following cases.

\textbf{Case 1: $\gamma(U_i)<\infty$ for all $0 \le i \le \gamma(u)$.} Since $\gamma(\N)=\infty$, it must be the case that $\gamma(U_i,U_j)=\infty$ for some $i\neq j$. Applying Lemma~\ref{canon13-biplem0} to the restriction of $\Delta$ to $G[U_i\cup U_j]$, we find a subset $X$ of $U_i\cup U_j$ such that $\gamma_G(X)=m$. 

\textbf{Case 2: $\gamma(U_i)=\infty$ for some $i\not=0$.} Let $H$ be the infinite subgraph of $G[U_i]$ obtained by deleting all the edges of colour $i$ from $G[U_i]$. Clearly, $\gamma_{H}(w)<\infty$ for all $w\in U_i$. So $H$ contains no infinite subset which is rainbow or star coloured. By the induction hypothesis, there is a subset $Y$ of $U_i$ such that $\gamma_{H}(Y)=m-1$. Observe that all the edges between $u$ and $Y\subset U_i$ have colour $i$, and since the colour $i$ is not counted by $\gamma_{H}$, we see that $\gamma_{G}(Y\cup\{u\})=m$. Therefore, $X=Y\cup\{u\}$ is the required subset.

\textbf{Case 3: $\gamma(U_0) = \infty$.} Let $V_0$ be the set of those vertices of $U_0$ not adjacent to $v$ in $G$. Since $\gamma(v)<\infty$, we have a partition of $U_0 \setminus V_0 = V_1\cup\dots\cup V_n$, with $n \le \gamma(v)$, based on the colour of the edge joining a given vertex of $U_0 \setminus V_0$ to the vertex $v$. Applying the same argument as in Cases 1 and 2 to the vertex $v$, we see that we are done unless $\gamma(V_0)=\infty$.  In this case, we consider the infinite subgraph $H$ of $G[V_0]$ obtained by deleting all the edges of colour $c$ from $G[V_0]$.

The fact that $\gamma_{G}(w)<\infty$ for all $w\in \N$ implies that $\gamma_{H}(w)<\infty$ for all $w\in V_0$. So $H$ has no infinite rainbow or star coloured subset. By the induction hypothesis, there is a subset $Y$ of $V_0$ such that $\gamma_{H}(Y)=m-1$. Observe that $uv$ has colour $c$ and there are no edges between $\{u,v\}$ and $Y\subset V_0\subset U_0$ in $G$. Since the colour $c$ is not counted by $\gamma_{H}$, it follows that $\gamma_{G}(Y\cup\{u,v\})=m$. Therefore, $X=Y\cup\{u,v\}$ is the required subset. This completes the proof.
\end{proof}

\section{Extensions and applications}
\label{canon13-extensions}

In this section, we shall first describe a finitary analogue of Theorem~\ref{canon13-canresult}. We then use this to prove Theorem~\ref{canon13-mainfvkm}. For us, a \emph{countable set} is a set that is either finite or countably infinite.

\subsection{Finitary extensions}

We can prove a version of Theorem~\ref{canon13-canresult} for colourings (of finite or infinite complete graphs) that use only finitely many colours.

\begin{theorem}\label{canon13-canresult'}
For all $n\in\N$, there exists a natural number $K = K(n)$ such that for every $k$-colouring $\Delta:V^{(2)}\twoheadrightarrow [k]$ of the complete graph on a countable set $V$ with $k \ge K$ colours, either
\begin{enumerate}
\item there is an $m$-coloured complete subgraph for every $m\in [n]$, or
\item there exists a rainbow coloured complete subgraph on $n$ vertices, or
\item there exists a star coloured complete subgraph on $n$ vertices.\qed
\end{enumerate}
\end{theorem}

This result can be proved by arguments similar to those used to prove Theorem~\ref{canon13-canresult}. There are two essential differences. First, as opposed to Theorem~\ref{canon13-erdos}, we use the following extension of the theorem proved by Erd\H{o}s and Rado, to colourings of \emph{finite} complete graphs with an arbitrary set of colours. 

\begin{theorem}\label{canon13-erdos'}
For every $n\in\N$, and every colouring $\Delta$ of the complete graph on a sufficiently large countable set $V$, there exists a subset $X$ of $V$ of size at least $n$ such that either
\begin{enumerate}
\item $X$ is $1$-coloured, or 
\item $X$ is rainbow coloured, or 
\item $X$ is lexically coloured.\qed
\end{enumerate}
\end{theorem}

Second, in the place of Lemma~\ref{canon13-biplem0}, we use the following finitary analogue which is proved in the same way as the lemma.
 
\begin{lemma}\label{canon13-biplem'}
For all $m,d\in\N$, there exists a natural number $L = L(m,d)$ with the following property: for every colouring $\Delta$ of a graph $G$ on a countable set $V$ such that
\begin{enumerate} 
\item $\gamma(v)<d$ for all $v\in V$, and 
\item there is a partition of $V=A\cup B$ such that $\gamma(A)<d$, $\gamma(B)<d$ and $\gamma(A,B)\ge L$, 
\end{enumerate}
there exists a subset $X$ of $V$ such that $X\cap A\not=\emptyset$, $X\cap B\not=\emptyset$ and $\gamma(X)=m$.\qed
\end{lemma}
 
\subsection{Applications}
 
Theorem~\ref{canon13-mainfvkm} may be deduced from Theorem~\ref{canon13-canresult'}. Recall that Theorem~\ref{canon13-mainfvkm} says for any natural number $n \in \N$, both $\binom{n}{2}, \binom{n+1}{2} \in \C{G}_\Delta$ for any colouring $\Delta$ of the complete graph on $\N$ using a finite, but sufficiently large number of colours.

\begin{proof}[Proof of Theorem~\ref{canon13-mainfvkm}]
We prove two propositions which, taken together, imply the result. The first is an easy corollary of Theorem~\ref{canon13-canresult'}

\begin{proposition}
\label{canon13-choose2}
For all $n\in\N$, there exists a natural number $C_1 = C_1(n)$ such that for any $k$-colouring $\Delta: V^{(2)} \twoheadrightarrow [k]$ of the complete graph on a countable set $V$ with $k \ge C_1$ colours, $\binom{n}{2} \in \C{G}_\Delta$.
\end{proposition}

\begin{proof} Take $C_1(n)=K(\binom{n}{2})$, where $K$ is as guaranteed by Theorem~\ref{canon13-canresult'}.
\end{proof}

The next proposition is perhaps not as straightforward.
\begin{proposition}
\label{canon13-1pchoose2}
For all $n\in\N$, there exists a natural number $C_2=C_2(n)$ with the property that for all $k \ge C_2$, there exists a natural number $D_{k,n}$ such that for any $k$-colouring $\Delta: V^{(2)} \twoheadrightarrow [k]$ of the complete graph on a countable set $V$ with $k \ge C_2$ colours, $\binom{n}{2}+1 \in \C{G}_\Delta$, provided $|V| \ge D_{k,n}$.
\end{proposition}

\begin{proof}
For $n=2$, it is an easy exercise to check that the result is true with $C_2(2)=2$ and $D_{k,2}=R(k+1;k)$, where $R(k+1;k)$ is the Ramsey number for finding a $1$-coloured copy of a complete graph on $k+1$ vertices when using $k$ colours.

For $n\ge 3$, let $s=n^4$. We claim that $C_2(n)=K(s)$ will do, where $K$ is the constant guaranteed by Theorem~\ref{canon13-canresult'}. For $k\ge C_2(n)$, we take $D_{k,n}=k^s+s+1$. Now, suppose that $\Delta:V^{(2)}\twoheadrightarrow[k]$ is a $k$-colouring and $|V|\ge D_{k,n}$. Then, by our choice of $C_2(n)$, either
\begin{enumerate}
\item there is an $m$-coloured complete subgraph for every $m\in [s]$, or
\item there exists a rainbow coloured complete subgraph on $s$ vertices, or
\item there exists a star coloured complete subgraph on $s$ vertices.
\end{enumerate}
Note that a star coloured complete subgraph on $s$ vertices contains an $m$-coloured complete subgraph for $2<m\le s$. Since $2<\binom{n}{2}+1\le s$, we are done unless there exists a rainbow coloured complete subgraph on $s$ vertices. Hence, suppose that the complete subgraph on the vertex set $S=\{u_1,u_2,\dots,u_s\}$ is rainbow coloured. For each $x\in V\setminus S$, there are $k^s$ possible values for the $s$-tuple $(\Delta(xu_1),\Delta(xu_2),\dots,\Delta(xu_s))$. Since, $|V \setminus S|\ge D_{k,n}-s>k^s$, we can find vertices $x, y\in V \setminus S$ such that 
\[
(\Delta(xu_1),\Delta(xu_2),\dots,\Delta(xu_s))=(\Delta(yu_1),\Delta(yu_2),\dots,\Delta(yu_s)).
\]

We claim that there is a subset $T\subset S$ of size $t=n^2$ such that for all $u\in T$, $\Delta(xu)\not\in \Delta(T^{(2)})$. Assume for the sake of contradiction that for every subset $T\subset S$ of size $t$, there exists at least one vertex $u\in T$ such that $\Delta(xu)\in \Delta(T^{(2)})$. Consider the set
\[
A=\{(u,T):u\in T\subset S,\,|T|=t,\,\Delta(xu)\in \Delta(T^{(2)})\}.
\]
By our assumption, for each $T\subset S$ of size $t$, there is at least one $u\in T$ such that $(u,T)\in A$, so $|A|\ge \binom{s}{t}$. As $S$ is rainbow coloured, there is at most one edge $ab$ in $S^{(2)}$ of colour $\Delta(xu)$ for each $u\in S$. If $(u,T)$ is in $A$, then we must have $a,b\in T$. So for each $u\in S$, there are at most $\binom{s-2}{t-2}$ sets $T$ such that $(u,T) \in A$. Thus, $|A|\le s\binom{s-2}{t-2}$. Combining these two inequalities for $|A|$, we get
\[
\binom{s}{t}\le |A| \le s\binom{s-2}{t-2}.
\]
This means that $t(t-1)\ge s-1$, contradicting the fact that $s=t^2$.

Hence, there is indeed a subset $T$ of $S$ of size $t=n^2$ such that $\Delta(xu)\not\in \Delta(T^{(2)})$ for all $u\in T$. Let $\C{Q}=\{\Delta(xu):u\in T\}$. If $|\C{Q}|<n$, then as $|T|=n^2$, there are vertices $v_1,v_2,\dots,v_n$ in $T$ such that
\[
\Delta(xv_1)=\Delta(xv_2)=\dots=\Delta(xv_n).
\]
Since this colour $\Delta(xv_1)$ is not an element of $\Delta(T^{(2)})$, we conclude that the set $\{x,v_1,v_2,\dots,v_n\}$ is $(\binom{n}{2}+1)$-coloured.

So we may assume that $|\C{Q}|\ge n$. Then there is a subset $U\subset T$ of size $n$ such that the colours $\Delta(xu)$ are distinct for all $u\in U$. Since $ U\subset T$, the colour $\Delta(xu)$ is not an element of $\Delta(U^{(2)})$ for each $u\in U$. We hence conclude that $U\cup\{x\}$ is rainbow coloured. 

Recall that there is a vertex $y\not=x$ in $V\setminus S$ such that $\Delta(xu)=\Delta(yu)$ for all $u\in S$. Since at most one edge $e$ in $(U\cup\{x\})^{(2)}$ is coloured with the same colour as the edge $xy$, by removing the endpoint of $e$ which lies in $U$ if necessary, we can find a subset $U'$ of  $U$ of size $n-1$ such that $\Delta(xy)$ is not an element of $\Delta((U'\cup\{x\})^{(2)})$. Then $U'\cup\{x,y\}$ is $(\binom{n}{2}+1)$-coloured since $U'\cup\{x\}$ and $U'\cup\{y\}$ are rainbow coloured sets of size $n$ using the same set of colours.
\end{proof}

It is easy to see that, taken together, Corollary~\ref{canon13-choose2} and Theorem~\ref{canon13-1pchoose2} imply Theorem~\ref{canon13-mainfvkm}.
\end{proof}

The following corollary of Lemma~\ref{canon13-biplem'} about finding $m$-coloured complete bipartite subgraphs might be of independent interest.

\begin{corollary}
For all $m\in\N$, there exists a natural number $B=B(m)$ such that if $\Delta:U\times V\twoheadrightarrow [k]$ is a $k$-colouring of the complete bipartite graph between two countable sets $U$ and $V$ with $k \ge B$ colours, then there exist $X\subset U$ and $Y\subset V$ such that the complete bipartite subgraph between by $X$ and $Y$ is $m$-coloured.
\end{corollary}

\begin{proof}
It is easy to verify that it suffices to to take $B(m)=L(m,m)$, where $L$ is the constant guaranteed by Lemma~\ref{canon13-biplem'}.
\end{proof}

\section{Conclusion}
\label{canon13-conclude}

We conclude by mentioning two questions that would merit further study. First, the problem of determining for each $k \in \N$, which natural numbers $m \in \N$ are guaranteed to belong to $\C{G}_\Delta$ for every $k$-colouring $\Delta:\N^{(2)} \twoheadrightarrow [k]$ is quite interesting; while we have taken a few steps towards this in this paper, the full question is still far from being resolved. Second, it would be reasonable to ask the questions considered here for $r$-uniform hypergraphs. However, even in the case of $\N^{(3)}$, it is not immediately clear  to us what the canonical structures analogous to the rainbow coloured and star coloured complete graphs should be.

\section*{Acknowledgements} The research in this paper was carried out while the authors were visitors at Microsoft Research, Redmond. We are grateful to Yuval Peres and the other members of the Theory Group at Microsoft Research for their hospitality. 

\bibliographystyle{amsplain}
\bibliography{canon_col}

\end{document}